\documentclass[sn-mathphys]{sn-jnl}

\theoremstyle{thmstyleone}%
\newtheorem{theorem}{Theorem}
\newtheorem{lemma}[theorem]{Lemma}%

\theoremstyle{thmstyletwo}%

\theoremstyle{thmstylethree}%

\newcommand{\Qq}{\mathbb{Q}}

\newcommand{\Zz}{\mathbb{Z}}

\begin{document}

\title[Subdividing triangles with $\pi$-commensurable angles]{Subdividing triangles with $\pi$-commensurable angles}


\author[1]{\fnm{Hasan} \sur{Korkmaz}}\email{hkorkmaz1@yahoo.com}

\author*[2]{\fnm{Ferit} \sur{\"{O}zt\"{u}rk}}\email{ferit.ozturk@boun.edu.tr}
\equalcont{These authors contributed equally to this work.}

\affil[1]{\orgdiv{Retired}, \orgname{\.Izmir Science High School},  \orgaddress{\city{\.Izmir}, \country{Turkey}}}

\affil*[2]{\orgdiv{Department of Mathematics}, \orgname{Bo\u{g}azi\c{c}i \"{U}niversitesi}, \orgaddress{\street{Bebek}, \city{\.Istanbul}, \postcode{34342}, \country{Turkey}}}


\abstract{
A point in the interior of a planar triangle determines a subdivision into six subtriangles.
A triangle with angles commensurable with $\pi$ is called $\pi$-commensurable.
For such a triangle a subdivision where each of the subtriangles are $\pi$-commensurable too 
is called $\pi$-commensurable.
We prove that there are infinitely many $\pi$-commensurable triangles that do not admit any $\pi$-commensurable subdivision
except the one given by angle bisectors.  
We count the number of $\pi$-commensurable subdivisions of triangles.
We perform a similar count for $\Zz$-degree subdivisions  of $\Zz$-degree triangles too.
Finally we show that subdivision by angle bisectors is essential in recursive subdivisions in the sense that 
 recursive $\pi$-commensurable subdivisions of any $\pi$-commensurable triangle ultimately involve a subdivision by  angle bisectors.
}




\maketitle

\section{The problem}
Consider a triangle in the Euclidean plane with integer angles in degrees and a point in its interior. Connecting that point to the vertices of the triangle, six interior angles are formed. 
In such a setup the school teachers naturally wish to have these six angles integer valued 
in degrees too, if they want to produce  `reasonable' looking triangle questions. 
But is it true that every triangle with integer angles contains such a point in its interior?

\begin{figure}[h!]
\begin{center}
\includegraphics[width=12cm]{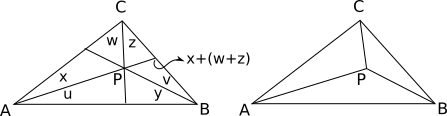}
\caption{Subdivisions determined by an interior point $P$}
\label{bolumle}
\end{center}
\end{figure}

Let us make this question slightly more general. 
A triangle in the Euclidean plane with its angles commensurable with $\pi$ is called a {\em $\pi$-commensurable triangle}.
We divide such a triangle into six  subtriangles  using an interior point as above.
If the six subtriangles (or equivalently three, as in the right figure) are $\pi$-commensurable too, we will call such a subdivision $\pi$-commensurable. 
If all the angles involved have integer degrees we will call the triangle and the subdivision {\em $\Zz$-degree}.
We ask how many $\pi$-commensurable (respectively $\Zz$-degree) subdivisions a $\pi$-commensurable  (respectively $\Zz$-degree) triangle has and whether there is any at all for a given triangle. It is known since the ancient times that the angle bisectors determine a $\pi$-commensurable subdivision for a $\pi$-commensurable  triangle.
In this short note we will show that except for the isosceles triangles, a triangle can have at most finitely many $\pi$-commensurable subdivisions and infinitely many have none
except via  angle bisectors. Furthermore most of the  $\Zz$-degree triangles do  not have a $\Zz$-degree subdivision.
We will give a complete list of those which admit no $\Zz$-degree subdivision and determine the number of distinct $\Zz$-degree subdivisions.

Now let $\tau$ be a triangle and  $P$ a point in its interior yielding a subdivision into 3 subtriangles. 
We set: $u=PAB/\pi$, $x=PAC/\pi$, $v=PBC/\pi$, $y=PBA/\pi$ and $w=PCA/\pi$, $z=PCB/\pi$; here the angles are in radians.
By the trigonometric Ceva's theorem it follows that these real numbers satisfy
\begin{eqnarray}
\label{ceva} \sin (u\pi) \sin (v\pi) \sin (w\pi) = \sin (x\pi) \sin (y\pi) \sin (z\pi), \\  \nonumber \mbox{with }u+v+w+x+y+z=1. 
\end{eqnarray}
Conversely if $u,v,w,x,y,z$ are positive real numbers that satisfy Equations~(\ref{ceva}) then they determine a unique triangle $\tau$. In fact, consider the unit circle  partitioned by six points $A_1,  B_1, C_1, A_2, B_2, C_2$ ordered positively so that they define six arcs measuring $2\pi$ times  $u,v,w,x,y,z$ in the same order. Then the triangle $A_1, B_1, C_1$ is the required one provided that the three diagonals $A_1B_2,A_2C_1,B_1C_2$ of the hexagon meet at a point (see \cite{PR}, page~138). But this is exactly the problem dealt with in \cite{PR}:  counting the number of multiple intersections of the diagonals of regular  polygons (or more generally those with interior angles $\pi$-commensurable).
To put it in other words, the relation of their problem with ours is via the observation that a $\pi$-commensurable 
triangle can be considered as a degenerate  $\pi$-commensurable hexagon in which every other interior angle measures $\pi$.

Having said these,  the number of distinct pairs $(\tau,P)$ with 
$\tau$ $\pi$-commensurable and the interior point $P$  yielding a $\pi$-commensurable subdivision  is equal to the number of distinct solutions 
of the ordered positive rational  numbers that satisfy  (\ref{ceva}). This latter problem has been termed as solving a trigonometric Diophantine equation in \cite{CJ} 
where it was proven that solving such an equation reduces to solving vanishing sums of roots of unity. In our problem, replacing the sines with complex exponentials we obtain 
a vanishing sum of 6 unit complex numbers and their conjugates. The requirement of rationality turns this sum into the following vanishing sum of roots of unity:
$$ \sum_{j=1}^{6} e^{i\pi\alpha_j} + \sum_{j=1}^{6} e^{-i\pi\alpha_j} = 0$$
with the extra condition that $\sum_{j=1}^{6} \alpha_j =1$ (see \cite[Equation~2.3]{PR}). 
In the sense of \cite{CJ} this is a vanishing sum of roots of unity of length 12.

Maybe the simplest instance of a trigonometric Diophantine equation is $\sin x = q\in\Qq$, where $x$ is to be $\pi$-commensurable.
The basic fact that there are exactly 3 solutions for the pair $(x,q)$ with $x\in[0,\pi/2]$ is referred as Niven's theorem (see \cite{Ni}, \cite{Le} or \cite{Ol}), 
although this result was probably known before  twentieth century. 
Trigonometric Diophantine equations giving way to particular vanishing sums of roots of unity with higher lengths have been worked out since then;
to name a few: \cite{Ne},  \cite{CJ} (lengths $\leq 9$), \cite{My} (particular cases with short lengths)  \cite{PR} (a particular case with length 12), 
\cite{LL} (a general result on lengths).
In particular the articles \cite{CJ} and \cite{PR} here  present, among other things, particular Diophantine equations that arise from intriguing geometric problems.

The complete list of solutions for (\ref{ceva}) has been given in \cite{PR} by explicitly counting the 
corresponding vanishing sums of roots of unity.  This result in turn gives an answer to our initial question.
In the next section we recall the results mentioned and give explicit counts for our problem. Furthermore we prove 
in Section~\ref{recurse} that 
one is ultimately forced to subdivide by  angle bisectors in a process of 
recursive $\pi$-commensurable subdivisions of any given $\pi$-commensurable triangle.

\section{The solution}
Let us summarize the solutions described in \cite[Theorem~4.4]{PR}. Given the triangle $\tau$ and its $\pi$-commensurable angles $A,B,C$ the first obvious family of solutions is any $\pi$-commensurable triple 
$u,v,w$ with sum $\pi/2$ and $x,y,z$ are a permutation of $u,v,w$. 
These solutions correspond to 6 pairs of square roots of unity adjusted to satisfy the required relation. Each of these solutions is coined as {\em trivial} in \cite{PR}
and denoted by $6R_2$. They are simply obtained  by solving the linear equation system 
\begin{equation*}
u+u'=A, \qquad v+v'=B, \qquad  w+w'=C,
\end{equation*}
where $(u',v',w')$ is any permutation of $(u,v,w)$. Clearly there are three types of solutions overall: either (i) $u=x=A/2$, $v=y=B/2$ and $w=z=C/2$; 
or (ii) without loss of generality $B=C$, i.e. $\tau$ is isosceles, and  $u=x=A/2, v+w=B, y=w, z=u$; 
or (iii) the sum of each pair among $A,B,C$ is greater than the one left out, i.e. $\tau$ is an acute triangle, which gives a unique solution for $u,v,w$. 
Solution (i) is exactly the subdivision determined by the angle bisectors of $\tau$. Thus any $\pi$-commensurable triangle has a trivial $\pi$-commensurable subdivision: via angle bisectors. In case of (ii) there are countably infinite subdivisions. Non-isosceles acute triangles have two more subdivisions via (iii).

As for the $\Zz$-degree setup, the cases above boil down to the following cases (angles in degrees):
either (i) $A,B,C$ are all even, or (ii) without loss of generality $A$ is even and $B=C$, or (iii) when $\tau$ is an acute triangle. These yield to (i) one, (ii) $B-1$, (iii) unique distinct trivial solutions for $u,v,w$.
As a conclusion, we have:
\begin{theorem}
\label{ikiz}
A $\Zz$-degree triangle has a $\Zz$-degree subdivision if and only if either it is isosceles or acute or its angles measure even degree.
\end{theorem}



The second type of solutions in  \cite{PR} is constituted of four (infinite) families for $(u,x,v,y,w,z)$ in degrees 
parametrized by $t\in\Qq$. To use in the sequel we list them below, in degrees:
\begin{subequations}\label{aile}
\begin{align}
(30,60+t,t,t,60-2t,30-t), 0<t<30 \label{aile1}\\
 (30,30-t,90-3t,2t,t,30+t), 0<t<30 \label{aile2} \\
 (30,30-2t,30-2t,t,2t,90+t), 0<t<15 \label{aile3} \\
 (60-4t,30-2t,t,3t,60+t,30+t), 0<t<15 \label{aile4}.
\end{align}
\end{subequations}

Obviously there are infinitely many $\pi$-commensurable triangles left out of these families.
Meanwhile there are (finitely) many integer solutions among these families.

The last type of solutions for (\ref{ceva}) are the 65 {\em sporadic} ones listed in \cite[Table~4.2]{PR}, which are not included in any of the previous families. The finiteness of these and the on-going discussion proves
\begin{theorem}
\label{piciler}
Any $\pi$-commensurable triangle has a $\pi$-commensurable subdivision, which is  determined by the angle bisectors. Moreover there are infinitely many $\pi$-commensurable triangles which do not have any other $\pi$-commensurable subdivision. Except for isosceles triangles, every $\pi$-commensurable triangle has finitely many $\pi$-commensurable subdivisions.
\end{theorem}

We note the following observation for future use:
\begin{lemma}
\label{genis}
The largest possible angle in the triangles obtained via families in (\ref{aile}) is 135. 
The largest and smallest possible angles in the sporadic triangles are 150 and 180/21 respectively. (Angles in degrees.)
\end{lemma}

Going through all possible $\Zz$-degree triangles that families in~(\ref{aile}) and in  \cite[Table~4.2]{PR} yield,  
we observe that there are cases  which have not been covered by Theorem~\ref{ikiz}, i.e. non-isosceles, non-acute triangles with at least one angle odd and which admit $\Zz$-degree subdivisions.  Performing a  calculation on computer we see that there are 119 such triangles up to symmetry.
Meanwhile the number of non-isosceles, non-acute triangles with at least one angle odd can be computed to be 1496 up to symmetry. This shows
\begin{theorem}
There are 1377 distinct $\Zz$-degree triangles up to similarity which do not admit any $\Zz$-degree
subdivision. All of these are non-isosceles, non-acute triangles with at least one angle of odd degree.\footnote{The list of  $\Zz$-degree triangles that do not admit any $\Zz$-degree subdivision can be found in \\ \url{www.math.boun.edu.tr/instructors/ozturk/noZDegreeSubdivision.html}.
}
\end{theorem}

Note that there are 2700 $\Zz$-degree triangles up to similarity. 

\subsection{Recursive subdivisions}
\label{recurse}
If one wishes to recursively subdivide a triangle in a $\pi$-commensurable fashion, it turns out that the consistent way of doing that is via angle bisectors. By saying this we mean the following observation:

\begin{theorem}
Recursive $\pi$-commensurable subdivisions of a $\pi$-commensurable triangle after sufficiently many steps produce at least one subtriangle which does not assume a $\pi$-commensurable subdivision except via angle bisectors.
 \end{theorem}

To prove this, first note that after sufficiently many subdivisions each subtriangle sharing a vertex with the original triangle has an angle arbitrarily small.
Let us call a triangle marginal if it is both obtuse and has an angle measuring less than, say, 1 in degrees.
We will observe that any subdivision of any subtriangle with a small angle ($<1$) produces a marginal subtriangle. 
By Lemma~\ref{genis}, a marginal triangle does not have a $\pi$-commensurable subdivision except the trivial subdivisons.
Moreover since a marginal triangle is neither acute nor isosceles, its only possible subdivison is via the angle bisectors and the proof follows.

Thus consider a subtriangle with a small angle ($<1$). If it has a nontrivial subdivision, it follows from Lemma~\ref{genis} that the triangle must belong to the families in (\ref{aile}). For example consider the triangle that appears in family~(\ref{aile1}) with its angle $v+y=2t$ close to 0. Then $t$ must be small.
The triangle has either the angles $u+x,v+y,w+z$ or $u+z,v+y,w+x$. In the former case, after one further subdivision suggested by (\ref{aile1}) we obtain a subtriangle with two angles $x+(w+z)=150-2t$ and $v=t$. (See Figure~\ref{bolumle}. Here the parantheses show that $w+z$ is an angle of the initial triangle.) 
This is exactly a marginal subtriangle for small $t$.
Similarly in the latter case, after one further subdivision offered by (\ref{aile1}) we obtain a subtriangle with two angles $(x+w)+z=150-2t$ and $v=t$. Going through all the other triangles with small angles,  we detect the following subtriangles:
\begin{eqnarray*}
\label{kucukacili} 
\nonumber \mbox{from }(\ref{aile1}): & w+z=90-3t, t\sim 30; & y+(u+x)=90+2t\sim 150, z=30-t\sim 0 \\
\nonumber & w+z=90-3t, t\sim 30; & (y+u)+x=90+2t\sim 150, w=60-2t\sim 0 \\
\nonumber \mbox{from }(\ref{aile2}): & y+w=3t, t\sim 0; & (v+z)+u=150-2t\sim 150, y=2t\sim 0 \\
\nonumber & y+w=3t, t\sim 0; & (v+x)+u=150-4t\sim 150, y=2t\sim 0 \\
\nonumber \mbox{from }(\ref{aile2}): & x+v=120-4t, t\sim 30; & (z+w)+y=30+4t\sim 150, v=90-3t\sim 0 \\
\nonumber & x+v=120-4t, t\sim 30; & z+(w+y)=30+4t\sim 150, v=90-3t\sim 0 \\
\nonumber \mbox{from }(\ref{aile3}): & y+w=3t, t\sim 0; & z+(u+x)=150-t\sim 150, w=2t\sim 0 \\
\nonumber & y+w=3t, t\sim 0; & (z+u)+x=150-5t\sim 150, w=2t\sim 0 \\
\nonumber \mbox{from }(\ref{aile3}): & x+v=60-4t, t\sim 15; & (z+w)+y=90+4t\sim 150, v=30-2t\sim 0 \\
\nonumber & x+v=60-4t, t\sim 15; & z+(w+y)=90+4t\sim 150, v=30-2t\sim 0 \\
\nonumber \mbox{from }(\ref{aile4}): & v+y=4t, t\sim 0; & w+ (z+u)=150-2t\sim 150, y=3t\sim 0 \\
\nonumber & v+y=4t, t\sim 0; &  (w+z)+x=120, w=2t\sim 0 \\
\nonumber \mbox{from }(\ref{aile4}): & u+x=90-6t, t\sim 15; &  (v+y)+w=60+5t\sim 135, u=60-4t\sim 0 \\
\nonumber & u+x=90-6t, t\sim 15; &  v+(y+w)=60+5t\sim 135, x=30-2t\sim 0 
\end{eqnarray*}
All these detected subtriangles are marginal except the last three; The last three do not have any sporadic subdivision. The subtriangle third from the bottom can be immediately  checked  having no subdivision coming from the families in (\ref{aile}) either. Finally, the last two subtriangles do have subdivisons in the family (\ref{aile3}). 
However those subdivisions produce at least one marginal subtriangle as the running discussion shows. Hence the proof.

Recall that there is an algebraic structure underneath the solutions of our problem \cite{PR}. 
It would be interesting to give algebraic reasons for the facts that we have brought forth.

\subsection{Synthetic proofs}
Let us go back to where we started off.  Consider a $\Zz$-degree subdivision of a triangle. Given several of  the subangles $u,v,w,x,y,z\in\Zz$ in degrees can you find the  measure of the remaining ones in a synthetic way? Note that using the list of all $\Zz$-degree subdivisions of all $\Zz$-degree triangles, we may ask thousands of such questions.\footnote{A list can be found in \\ \url{www.hasankorkmaz-ifl.com/dosyalar/yaptigim-calismalar/ZSubdivisions.pdf}.
} 
We invite the readers to construct this sort of questions and attempt to solve them synthetically.

\end{document}